\documentclass[11pt]{amsart}
\usepackage{graphicx,calc,amscd}
\usepackage{epsfig,psfig}
\usepackage{float}
\usepackage{labelfig}
\usepackage{amsmath}
\usepackage{amsfonts}
\usepackage{amssymb}

\newcommand{\R}{{\mathbb R}}

\newcommand{\N}{{\mathbb N}}

\newcommand{\Z}{{\mathbb Z}}

\newcommand{\GGS}{{\mathbb G}(S)}
\newcommand{\CC}{{\mathcal C}}

\newcommand{\Sp}{{\mathbb S}}

\newcommand{\ii}{\mathrm{int}}

\newcommand{\arcsinh}{\mathrm{arcsinh}}
\newcommand{\arccosh}{\mathrm{arccosh}}

\newtheorem{theorem}{Theorem}[section]
\newtheorem{lemma}[theorem]{Lemma}
\newtheorem{corollary}[theorem]{Corollary}
\newtheorem{proposition}[theorem]{Proposition}

\newtheorem{definition}[theorem]{Definition}
\begingroup\makeatletter\ifx\SetFigFont\undefined%
\gdef\SetFigFont#1#2#3#4#5{%
 \reset@font\fontsize{#1}{#2pt}%
 \fontfamily{#3}\fontseries{#4}\fontshape{#5}%
 \selectfont}%
\fi\makeatother\endgroup%
\begin{document}
\title{Minimal length of two intersecting simple closed geodesics}
\author{Thomas Gauglhofer and Hugo Parlier}
\address{
   Section de Math\'ematiques, EPFL, 1015 Lausanne, SWITZERLAND}
\email{thomas.gauglhofer@a3.epfl.ch }

\address{Section de Math\'ematiques\\
Universit\'e de Gen\`eve\\
1211 Gen\`eve\\
SWITZERLAND} \email{hugo.parlier@math.unige.ch}

\thanks{The first author was supported in part by SNFS grant number
2100-065270, the second author was supported by SNFS grant number
PBEL2-106180.} 
\date{\today}
\begin{abstract}On a hyperbolic Riemann surface, given two simple closed geodesics
that intersect $n$ times, we address the question of a sharp lower
bound $L_n$ on the length attained by the longest of the two
geodesics. We show the existence of a surface $S_n$ on which there
exists two simple closed geodesics of length $L_n$ intersecting
$n$ times and explicitly find $L_n$ for $n\leq 3$.
\end{abstract}
\maketitle

\section{Introduction}                  \label{Sect:S1}

Extremal hyperbolic Riemann surfaces for a variety of geometric
quantities are objects of active research. Well known cases
include the study of surfaces with maximum size systoles
(\cite{ak03}, \cite{bav97}, \cite{sc931}), surfaces with largest
embedded disk (\cite{bav96}, \cite{gigo021}) or more classically
surfaces with maximum number of automorphisms (the study of
Hurwitz surfaces and related topics). These subjects are related
to the study of the simple length spectrum of a surface $S$,
denoted $\Delta_0(S)$, which is the ordered set of lengths of
(non-oriented, primitive) simple closed geodesics (with
multiplicity). The question of interpreting the geometry of the
surface through the values found in the simple length spectrum
seems to be a very difficult subject. For instance, for surfaces
of genus $2$ and $3$, it is not known whether the simple length
spectrum determines the surface up to isometry (see
\cite{bubook}).\\

One of the major tools used to approach these problems is the
collar theorem, and in particular a corollary which states that
two short simple closed geodesics (of length less than
$2\,\arcsinh(1)$) cannot intersect (see \cite{bubook}). The bound
is sharp because it can be realized on a particular torus with a
cusp. The bound is never reached for a closed surface, but for any
genus, is realized in
the compactification of its Moduli space.\\

The goal of this article is to generalize this result by studying
the relationship between the number of intersection points between
two simple closed geodesics and the length of the geodesics. The
surfaces we consider lie in the Moduli space of surfaces with
boundary $\mathcal{M}_{g,k}$, where $g$ is the genus and $k$ is
the number of simple closed boundary geodesics which we allow to
be cusps (geodesics of length $0$). The foundation of our study is
found in the following theorem (section 2).

\begin{theorem}\label{thm:base} On a hyperbolic Riemann surface $S$, let $\alpha$ and $\beta$
be simple closed geodesics that intersect $n$ times. Then there
exists a universal constant $L_n$ such that
$\max\{\ell(\alpha),\ell(\beta)\}\geq L_n$ and $L_n\longrightarrow
\infty$ when $n\longrightarrow \infty$. Furthermore, a surface
$S_n$ realizing the bound exists.
\end{theorem}

By realizing the bound, we mean that on $S_n$ there are two simple
closed geodesics of length $L_n$ that intersect $n$ times. We
further investigate the asymptotic behavior of $L_n$ in the
following proposition where we prove:

\begin{proposition}\label{prop:asymptoticintro}
Let $l_n$ be the positive solution of the equation
$$l_n=2n\,\arcsinh\left(\frac{1}{\sinh(l_n/2)}\right).$$

Then
$$l_n\leq L_n<2l_n.$$

\end{proposition}

We are able to describe the surfaces explicitly for $n\in
\{2,3\}$, which gives us the following result.

\begin{theorem} The surfaces $S_2$ and $S_3$ are once-punctured tori and $$\textstyle L_2 = 2\,\arccosh(2),\qquad L_3 = 2\,\arccosh
\left(\sqrt{\frac12\left(7+\frac{11}{3}\sqrt{\frac{11}{3}}\right)}\;\right).$$
\end{theorem}

As mentioned earlier, $L_1=2\,\arcsinh 1$ and note that the value
for $L_2$ was previously proved in \cite{gase05} but the proof
presented here is new. The surfaces $S_1$, $S_2$ and $S_3$ are all
distinct once-punctured tori. We show that they all have
non-trivial isometry groups. (A once-punctured torus is
necessarily hyperelliptic, so by non-trivial isometry group we
mean an isometry group not isomorphic to $\Z_2$). It seems
reasonable to conjecture that for all $n$, $S_n$ is also a
once-punctured torus with a non-trivial isometry group.\\

This article is organized as follows. Section \ref{sect:pre} is
devoted to preliminaries and the proof of theorem \ref{thm:base}
and proposition \ref{prop:asymptoticintro}. The next two sections
concern the exact values of $L_2$ and $L_3$ and are similar in
nature. The final section discusses possible future directions for
related questions.

\section{Preliminaries and groundwork}\label{sect:pre}

We will be considering hyperbolic Riemann surfaces of finite area,
with or without boundary. We allow boundary to be either cusps or
simple closed geodesics. A {\it surface} will always designate a
surface of this type. The signature of a surface will be denoted
$(g,k)$ where $g$ is the genus and $k$ the number of boundary
geodesics (or cusps). A surface of signature $(0,3)$ is called a
pair of pants, a surface of signature $(1,1)$ a one-holed torus,
and a surface of signature $(0,4)$ a four-holed sphere. We reserve
the term {\it punctures} for cusps, and holes can be cusps as well
as boundary geodesics. The Moduli space of surfaces with boundary
will be denoted $\mathcal{M}_{g,k}$.

The set of interior simple closed geodesics of a surface $S$ will
be denoted $\GGS$. The length of a simple closed geodesic $\alpha$
will be denoted $\ell(\alpha)$, although a geodesic and its length
might not be distinguished. The term {\it geodesic} will sometimes
be used instead of simple closed geodesic, but only if it is
obvious in the context. Closed geodesics will be considered to be
non-oriented (unless specified) and primitive, meaning that a
closed geodesic cannot be written as the $k$-fold iterate of
another closed geodesic. Seen this way, geodesics are point sets
independent of parameterization. We denote $\ii(\alpha,\beta)$ the
number of transversal intersection points between two simple
closed geodesics $\alpha$ and $\beta$. We define the simple length
spectrum $\Delta_0(S)$ as the ordered set of lengths of all
interior simple closed geodesics. Notice that our definition takes
into account multiplicity, namely that if there are $n$ distinct
simple closed geodesics of $S$ with equal length $\ell$, then the
value $\ell$ will appear $n$ times in $\Delta_0(S)$. Consider two
surfaces $S$ and $\tilde{S}$ with simple length spectra
$\Delta_0(S)=\{\ell_1\leq \ell_2 \leq \hdots \}$ and
$\Delta_0(\tilde{S})=\{\tilde{\ell}_1\leq \tilde{\ell}_2 \leq
\hdots \}$. The notation $\Delta_0(S)<\Delta_0(\tilde{S})$ is an
abbreviation for $\ell_i < \tilde{\ell}_i$ for all $i\in \N^*$.\\

In order to describe the pasting of a simple closed geodesic, one
generally uses {\it twist parameters}. The only use we will have
of twist parameters is to describe what we call {\it without
twist} or {\it zero-twist} and {\it half-twist}. Recall that a
pair of pants has three disjoint unique simple geodesic paths
between distinct boundary geodesics, called perpendiculars, which
decompose the pair of pants into two isometric hyperbolic
right-angled hexagons. If two pairs of pants are pasted along a
geodesic $\alpha$ such that the endpoints of the perpendiculars
coincide, then we refer to a pasting with {\it zero-twist} or {\it
without twist}. The terminology is slightly different for
one-holed tori. Consider a pair of pants with two boundary
geodesics $\alpha_1$ and $\alpha_2$ of equal length, and paste
$\alpha_1$ and $\alpha_2$ together in order to obtain a one-holed
torus. If the common perpendicular $a$ between $\alpha_1$ and
$\alpha_2$ has its endpoints that coincide, then we refer to a
pasting with {\it zero-twist} or {\it without twist}. If the
endpoints of $a$ are diametrically opposite on the geodesic
formally known as $\alpha_1$ or $\alpha_2$, then we refer to a
{\it half-twist}. Finally, if an interior simple closed geodesic
$\alpha$ is said to be {\it pasted with a half-twist}, then we
mean that it has been obtained from the construction described
above.\\

A function $f_\alpha:{\mathcal M}_{g,k} \longrightarrow \R^+$ that
associates to a closed geodesic $\alpha$ its length depending on
the choice of metric is generally referred to as a length
function. Length functions are well known to be analytic (one way
of seeing this is via Fricke trace calculus, see for instance
\cite{abbook}). What is interesting to us is that the function of
an interior closed geodesic remains continuous, even if boundary
length goes to $0$.\\

The collar theorem (i.e. \cite{kee74}, \cite{bu78}, \cite{ra79})
gives a very precise description of the\linebreak[4] geometry of
surfaces around simple closed geodesics.\pagebreak

\begin{theorem}\label{thm:collar} Let $\gamma_1$ and $\gamma_2$ be non-intersecting
simple closed geodesics on $S$. Then the collars
$$\CC(\gamma_i)=\{p\in S \mid d_S(p,\gamma_i)\leq w(\gamma_i)\}
$$
of widths
$$
w(\gamma_i)=\arcsinh(1/\sinh\frac{\ell(\gamma_i)}{2})
$$
are pairwise disjoint for $i=1,2$. Furthermore, each
$\CC(\gamma_i)$ is isometric to the cylinder
$[-w(\gamma_i),w(\gamma_i)]\times \Sp^1$ with the metric
$ds^2=d\rho^2+\ell(\gamma_i)^2\cosh^2\!\rho \, dt^2$.
\end{theorem}

This implies that a geodesic $\alpha$ that transversally
intersects $n$ times another geodesic $\beta$ satisfies
$\ell(\alpha)\geq 2n\,w(\beta)$.

\begin{corollary}\label{cor:L1} If $\alpha$ and $\beta$ are two simple
closed geodesics such that $$\ell(\alpha)\leq \ell(\beta)\leq
2\,\arcsinh 1,$$ then they do not intersect.
\end{corollary}

We shall also use the following result found in \cite{pa044}.

\begin{theorem}\label{thm:cusps}
Let $S$ be a surface of signature $(g,k)$ with $k > 0$. Let
$\gamma_1,\hdots,\gamma_k$ be the boundary geodesics of $S$. For
$(\varepsilon_1,\hdots,\varepsilon_k) \in (\R^+)^k$ with at least
one $\varepsilon_i\neq 0$, and $\varepsilon_j\leq \ell(\gamma_j)$
for all $j$, there exists a surface $\tilde{S}$ with boundary
geodesics of length
$\ell(\gamma_1)-\varepsilon_1,\hdots,\ell(\gamma_k)-\varepsilon_k$
such that all corresponding simple closed geodesics in $\tilde{S}$
are of length strictly less than those of $S$
$(\Delta_0(\tilde{S})<\Delta_0(S))$.
\end{theorem}

There is an immediate corollary to this result which is very
useful to our study.

\begin{corollary}\label{cor:cusp} If $\alpha$ and $\beta$ are two simple closed
geodesics that intersect $n$ times on a surface with non-empty
boundary and with at least one boundary geodesic not a cusp, then
there exists a surface of same signature, with only cusps as
boundary, containing two simple closed geodesics $\tilde{\alpha}$
and $\tilde{\beta}$ which intersect $n$ times and such that
$\ell(\tilde{\alpha})<\ell(\alpha)$ and
$\ell(\tilde{\beta})<\ell(\beta)$.
\end{corollary}

This corollary implies that we can limit ourselves to studying
surfaces with cusp boundary. This will be useful in the following
theorem which is the starting point of our work.

\begin{theorem}\label{thm:Ln} There exists a
universal constant $L_n$ such that
$\max\{\ell(\alpha),\ell(\beta)\}\geq L_n$ for any two simple
closed geodesics $\alpha$ and $\beta$ that intersect $n$ times on
a hyperbolic compact Riemann surface. Furthermore, a surface $S_n$
realizing the bound exists. Finally, $L_n\longrightarrow \infty$
when $n\longrightarrow \infty$.
\end{theorem}

\begin{proof}
The idea of the proof is to show that, for every $n$, we are
evaluating a continuous function on a finite set of compact sets.
The function is the one that associates to a surface $S$ the
following value:

$$
f(S)=\min_{\{\alpha,\beta \in \GGS\mid
\ii(\alpha,\beta)=n\}}\max\{\ell(\alpha),\ell(\beta)\}.
$$

For a given signature $(g,k)$, $f:{\mathcal
M}_{g,k}\longrightarrow \R^+$, is obviously continuous and
bounded. (Mind that for certain signatures, $f$ may not be
defined, for instance on surfaces of signature $(0,4)$, there are
no pairs of simple closed geodesics that intersect an odd number
of times.) Suppose $\alpha$ and $\beta$ are two simple closed
geodesics on a surface $S$ that intersect $n$ times. Consider the
following subsurface $S_{\alpha,\beta}$ of $S$. $S_{\alpha,\beta}$
is the surface (possibly with boundary) that embeds a tubular
neighborhood around the point set $\alpha\cup \beta$ such that all
interior simple closed geodesics of $S_{\alpha,\beta}$ intersect
either $\alpha$ or $\beta$. In other words, $\alpha$ and $\beta$
fill $S_{\alpha,\beta}$. For example, if $n=1$ this is necessarily
a surface of signature $(1,1)$. It is easy to see that the
signature $(g,k)$ of $S_{\alpha,\beta}$ is universally bounded by
a function of $n$ ($g\leq f_1(n)$, $k\leq f_2(n)$). There are thus
a finite number of possible signatures for $S_{\alpha,\beta}$,
which we shall denote $(g_1,k_1),\hdots (g_m,k_m)$. As any
interior simple closed geodesic of $S_{\alpha,\beta}$ intersects
either $\alpha$ or $\beta$, and as we are trying to minimize the
lengths, the collar theorem ensures that the length of the systole
of $S_{\alpha,\beta}$ is bounded from below (otherwise the maximum
length of $\alpha$ and $\beta$ would be unbounded). Denote by
$\epsilon_n$ this lower bound. By corollary \ref{cor:cusp}, as we
are searching for a minimal value among all surfaces, we can limit
ourselves to searching among surfaces with all boundary geodesics
being cusps. Denote by ${\mathcal M}_{g,k}^0$ the restricted set
of surfaces of signature $(g,k)$ with cusp boundary. Further
denote by ${\mathcal M}_{(g,k),\epsilon_n}^0$ the subset of
${\mathcal M}_{g,k}^0$ with systole bounded below by $\epsilon_n$.
We are now searching among a finite set of such sets, namely for
each $(g_j,k_j)$, $j\in \{1,\hdots,m\}$, we need to study the set
${\mathcal M}_{(g_j,k_j),\epsilon_n}^0$. These sets are well-known
to be compact (for surfaces without boundary see \cite{mu71}, and
with boundary see \cite{mat76}). As $f$ is a continuous function
that we allow to cover a finite number of compact sets, it follows
that $f$ admits a minimum, and the value of $f$ in this point we
denote
$L_n$. A point in Moduli space which reaches the minimum is denoted
$S_n$.\\

We now need to show that $L_n\longrightarrow \infty$ when
$n\longrightarrow \infty$. Suppose this is not the case, meaning
there exists some $L$ such that $L_n<L$ for all $n$. This would
mean that for any $n$, there exist two simple closed geodesics
$\alpha_n$ and $\beta_n$ on some surface $S$ that intersect $n$
times such that $\ell(\alpha_n)\leq \ell(\beta_n)\leq L$. By the
collar theorem $\ell(\beta_n)\geq 2n\,
\arcsinh(\frac{1}{\arcsinh(L/2)})$. But this is a contradiction,
because for any $L$, $n$ can be chosen so that this is not the
case. The theorem is now proven.
\end{proof}

The notation $S_{\alpha,\beta}$, as defined in the previous proof,
will be regularly used throughout the article. To study the
asymptotic behavior of $L_n$, we shall use the quantity $l_n$
defined in the following proposition.

\begin{proposition}\label{prop:l_n-increasing} For $n\in\mathbb{N}$, let $l_n$ be the
positive solution of the equation $$l_n=2n
\,\arcsinh\left(\frac{1}{\sinh(l_n/2)}\right).$$ Then $l_n$ is
strictly increasing in $n$.
\end{proposition}
\begin{proof}
The equation $l_n=2n \,\arcsinh(\frac{1}{\sinh(l_n/2)})$ is
equivalent to $\sinh(\frac{l_n}{2n})\sinh(\frac{l_n}{2})=1$.
Suppose now that there is an $n\in\mathbb{N}$ such that $l_n\geq
l_{n+1}$. Then $\sinh(\frac{l_n}{2})\geq\sinh(\frac{l_{n+1}}{2})$
which implies therefore that
$\sinh(\frac{l_n}{2n})\leq\sinh(\frac{l_{n+1}}{2n+2})$.

 But
$\frac{l_n}{2n}\leq\frac{l_{n+1}}{2n+2}$ implies $l_n<l_{n+1}$
which leads to a contradiction.
\end{proof}

Now the asymptotic behavior of $L_n$ can be expressed as follows.

\begin{proposition}\label{prop:l_n<L_n<2l_n}
$l_n\leq L_n<2l_n.$
\end{proposition}
\begin{proof}
Let us begin by showing $l_n\leq L_n$.

If a simple closed geodesic $\alpha$ of length $l_n$ intersects a
simple closed geodesic $\beta$ $n$ times, then $\beta$ is at least
as long as $2n$ times the width of the collar of $\alpha$. Thus
$\ell(\beta)\geq l_n$. The width of the collar of $\alpha$
increases when $\alpha$ gets shorter, thus $l_n\leq L_n$.\\

It remains to show that $L_n<2l_n$.

For $n\in\mathbb{N}$, let $\mathcal{Y}$ be a pair of pants whose
boundary consists of a cusp and two boundary geodesics, $\alpha_1$
and $\alpha_2$, both of length $l_n$. Let us paste these two
geodesics together (denote the resulting geodesic $\alpha$)
without twist. The common perpendicular between $\alpha_1$ and
$\alpha_2$ is now a simple closed geodesic, which we shall denote
$\delta$. Notice that $\ell(\delta)= l_n/n$. For a given primitive
parameterization of $\alpha$ and $\delta$, consider the simple
closed curve $\tilde{\beta}=\delta^n\alpha$ and its unique
geodesic representative $\beta$. By construction,
$\ell(\beta)<\ell(\tilde{\beta})=l_n+n(l_n/n)=2l_n$. We have thus
constructed a once-punctured torus with two interior geodesics
$\alpha$ and $\beta$ that satisfy $\ii(\alpha,\beta)=n$ and
$\max\{\ell(\alpha),\ell(\beta)\}<2l_n$. It follows that $L_n< 2
l_n$.
\end{proof}

Finally, as an illustration of our investigation, let us give the
value for $L_1$ and describe the surface $S_1$. Corollary
\ref{cor:L1} implies that $L_1\geq 2\,\arcsinh 1$. In
fact,\linebreak[4] $L_1=2\,\arcsinh 1$, and this can be shown by
constructing the surface $S_1$ which realizes the bound $L_1$.
Consider, in the hyperbolic plane, a quadrilateral with three
right angles and one zero angle (a point at infinity). This
quadrilateral can be chosen such that the two finite length
adjacent sides are of length $\arcsinh 1$. By taking four copies
of this quadrilateral, and pasting them together as in the
following figure, one obtains a once-punctured torus with two
simple closed geodesics of length $2\,\arcsinh 1$ that intersect
once.

\begin{figure}[h]
\leavevmode \SetLabels
\endSetLabels
\begin{center}
\AffixLabels{\centerline{\epsfig{file = 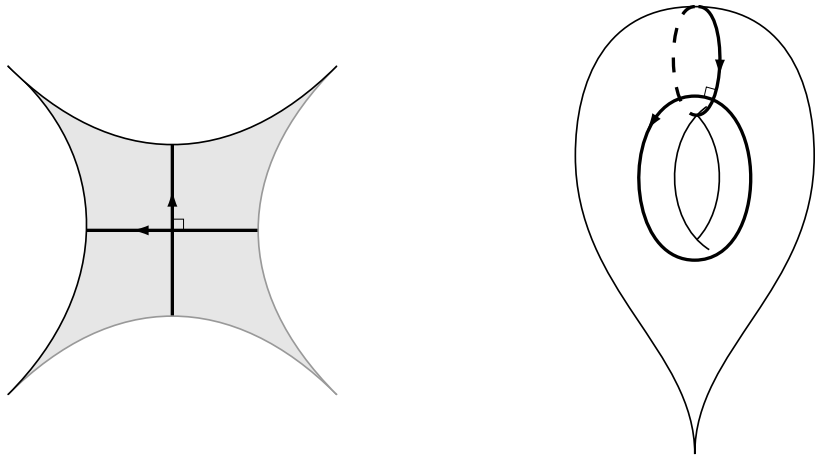,
height=6cm}}}
\end{center}
\caption{Surface $S_1$ with geodesics in bold}
\label{fig:surfaceS1}
\end{figure}

This once-punctured torus is the only surface on which two
intersecting geodesics can have length $L_1$. It is worth
mentioning that this torus has other remarkable properties: it is
the only once-punctured torus with an automorphism of order $4$.\\

\section{Finding $S_2$ and calculating $L_2$}

\begin{figure}[ht!]
\centerline{ \SetLabels
(.9*.2) {$\alpha$}\\
(.25*.2) {$\alpha$}\\
(.0*.53) {$\beta$}\\
(.545*.53) {$\beta$}\\
(.11*.4) {$A$}\\
(.26*.4) {$B$}\\
(.64*.4) {$A$}\\
(.805*.4) {$B$}\\
(.05*.95) {case 1}\\
(.5*.95) {case 2}\\
\endSetLabels
\AffixLabels{
 \hspace{.8mm}\psfig{figure=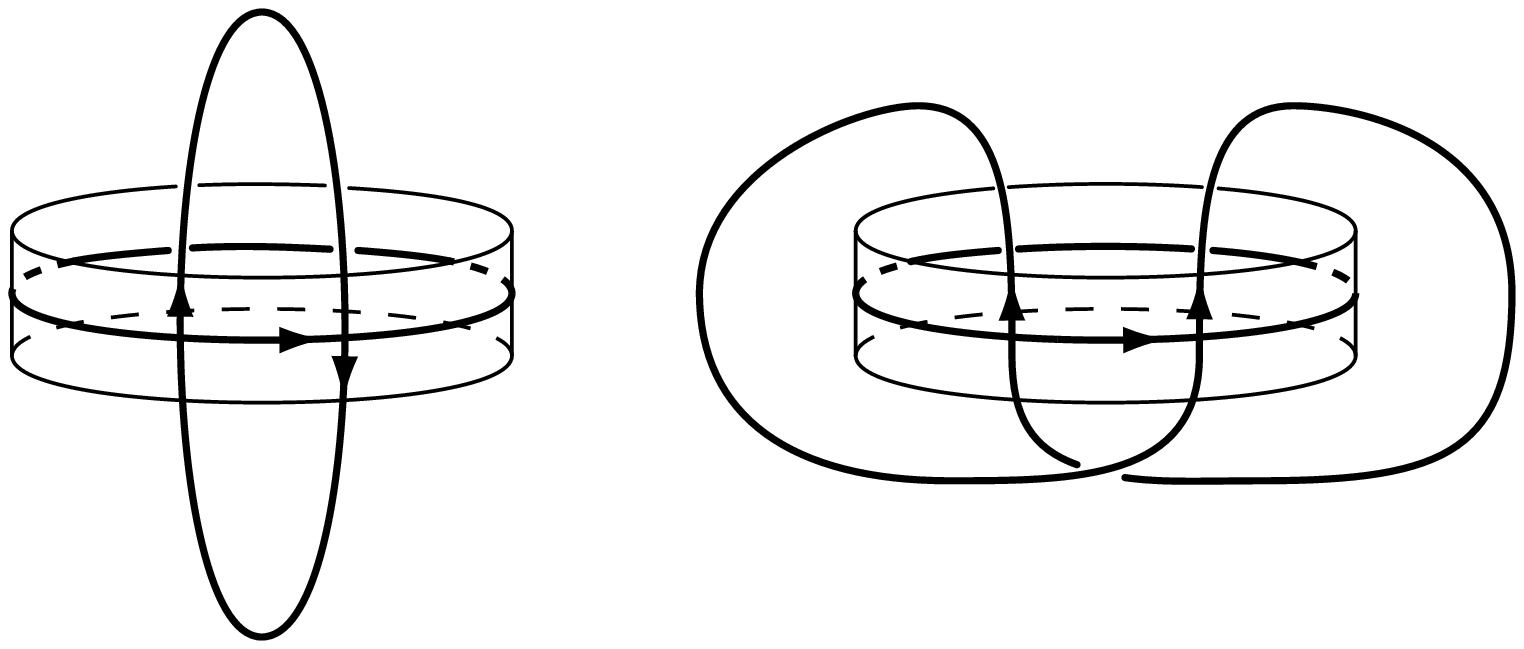,height=4.5cm}
}}\caption{The two cases for $n=2$} \label{fig:twoInter}
\end{figure}

Let us consider two simple closed geodesics $\alpha$ and $\beta$
on a surface $S$ that intersect twice in points $A$ and $B$, and
the subsurface $S_{\alpha,\beta}$ as defined in the proof of
theorem \ref{thm:Ln}. In order to distinguish possible signatures
for the surface $S_{\alpha,\beta}$, let us give $\alpha$ and
$\beta$ orientations. Let $C_\beta$ be a collar around $\beta$.
The ordered pair of simple closed \emph{oriented} geodesics
$(\alpha,\beta)$ induces an orientation on $C_\beta$ in both $A$
and $B$. These orientations are either opposite (case 1) or the
same (case 2). This is illustrated in figure \ref{fig:twoInter}.

\begin{figure}[ht!]
\centerline{ \SetLabels
(.35*-.1) {$\alpha$}\\
(.67*-.1) {$\alpha$}\\
(.63*-.1) {$\tilde\varepsilon_4$}\\
(.63*1.05) {$\tilde\varepsilon_3$}\\
(1*.56) {$\tilde\varepsilon_2$}\\
(0*.56) {$\tilde\varepsilon_2$}\\
(1*.46) {$\beta$}\\
(1*.36) {$\tilde\varepsilon_1$}\\
(0*.36) {$\tilde\varepsilon_1$}\\
\endSetLabels
\AffixLabels{
\hspace{2mm}\psfig{figure=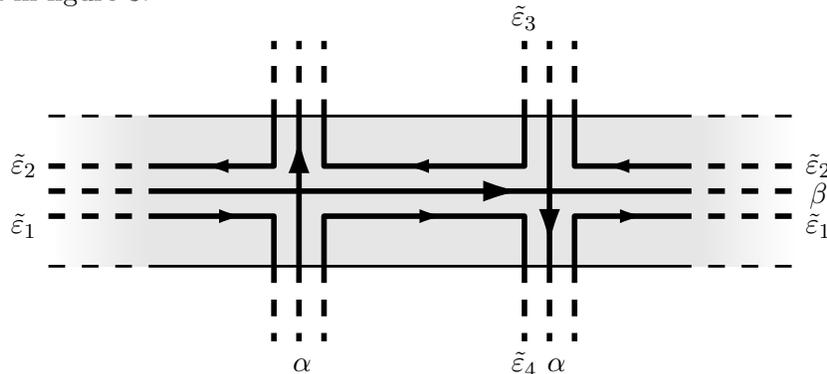,height=4cm}\hspace{2mm}}}
\caption{The simple closed curves $\tilde\varepsilon_1,
\tilde\varepsilon_2, \tilde\varepsilon_3$ and
$\tilde\varepsilon_4$ in case 1}\label{fig:AuBd}
\end{figure}

In case 1, $S_{\alpha,\beta}$ is a surface of signature $(0,4)$
obtained by cutting along the simple closed geodesics
$\varepsilon_1, \varepsilon_2, \varepsilon_3$ and $\varepsilon_4$
homotopic to the simple closed curves $\tilde\varepsilon_1,
\tilde\varepsilon_2, \tilde\varepsilon_3$ and
$\tilde\varepsilon_4$ shown in figure \ref{fig:AuBd}.

\begin{figure}[ht!]
\centerline{\SetLabels
(.29*.45) {?}\\
(.73*.45) {?}\\
(0.03*.11) {$\alpha$}\\
(.535*.11) {$\alpha$}\\
(.3*.99) {$\beta$}\\
(.41*.27) {$\tilde\varepsilon_1$}\\
(.6*.27) {$\tilde\varepsilon_2$}\\
(.05*-.07) {$A$}\\
(.05*.99) {$B$}\\
(.95*-.07) {$A$}\\
(.95*.99) {$B$}\\
(.5*.99) {$A$}\\
(.5*-.07) {$B$}\\
\endSetLabels
\AffixLabels{
\hspace{2mm}\psfig{figure=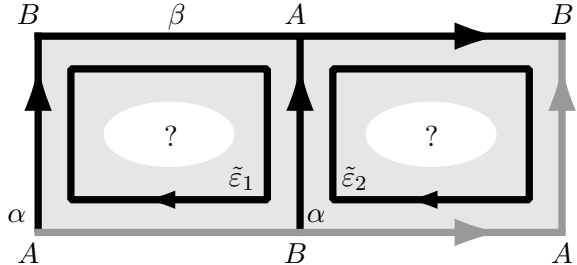,height=3cm}\hspace{2mm}}}
\caption{The simple closed curves $\tilde\varepsilon_1$ and
$\tilde\varepsilon_2$ in case 2}\label{fig:AuBu}
\end{figure}

In case 2, there are two possible topological situations for the
minimal surface $S_{\alpha,\beta}$. Indeed, consider the simple
closed curves $\tilde\varepsilon_1$ and $\tilde\varepsilon_2$
shown in Figure \ref{fig:AuBu}. One of these curves may be
null-homotopic, but not both because otherwise the surface would
be a torus without holes, which of course cannot admit a
hyperbolic metric.

If only one curve is not null-homotopic, say
$\tilde\varepsilon_1$,  we cut the surface along the geodesic that
is homotopic to $\tilde\varepsilon_1$ to obtain a surface of
signature $(1,1)$. If neither curve is null-homotopic, we cut the
surface along the two geodesics homotopic to $\tilde\varepsilon_1$
and $\tilde\varepsilon_2$ to obtain a surface of signature
$(1,2)$.\\

Therefore, in view of corollary \ref{cor:cusp}, $S_2$ is a sphere
with four cusps, a torus with one
cusp or a torus with two cusps.\\

First let us investigate geodesics intersecting twice on a
four-holed sphere.

\begin{proposition}\label{prop:fourholed2}
Let $X$ be a four-holed sphere (where we allow the boundary
geodesics to be cusps). Let $\alpha$ and $\beta$ be distinct
interior simple closed geodesics of $X$. Then
$$
\max\{\ell(\alpha),\ell(\beta)\}\geq 4\,\arcsinh 1 = 2\,\arccosh
3.
$$

Furthermore equality holds for a sphere with four cusps obtained
by gluing two pairs of pants with two cusps and third boundary
geodesic of length $2\,\arccosh 3$ without twist.
\end{proposition}

\begin{proof}
By corollary \ref{cor:cusp} it suffices to show the result for $X$
a sphere with four cusps. Suppose $\ell(\alpha)\geq\ell(\beta)$.
Now suppose by contradiction that $\ell(\alpha)<4\,\arcsinh 1$.
On a four-holed sphere, distinct interior simple closed geodesics
cross at least twice. By the collar theorem, the length of any
other interior simple closed geodesic must be strictly greater
than four times the width $w(\alpha)$ of the half-collar around
$\alpha$, which by \ref{thm:collar} is
$w(\alpha)=\arcsinh(\frac{1}{\sinh(2\,\arcsinh 1/2)})=\arcsinh 1$.
Thus $\ell(\beta)>4\,\arcsinh 1$, a contradiction. Thus equality
can only be attained if both $\alpha$ and $\beta$ are of length
$4\,\arcsinh 1$. It follows that a surface on which equality is
reached has a simple closed geodesic $\alpha$ of length
$4\,\arcsinh 1$. If there is any twist around this geodesic, then
all simple closed geodesics crossing $\alpha$ are of length
strictly superior to $4\,\arcsinh 1$ which concludes the argument.
\end{proof}

Let us now consider the case of two geodesics that intersect twice
on a one-holed torus. We recall that one-holed tori are
hyperelliptic, and we shall refer to the three interior fixed
points of the hyperelliptic involution as the Weierstrass points.

\begin{figure}[h]
\leavevmode \SetLabels
(.54*.78) $\alpha$\\
(.55*.32) $h_\alpha$\\
\endSetLabels
\begin{center}
\AffixLabels{\centerline{\epsfig{file=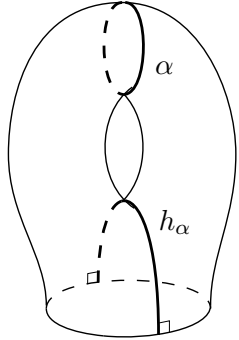,height=4.5cm}}}
\end{center}
\caption{The height $h_\alpha$ associated to $\alpha$}
\label{fig:qpieceheight}
\end{figure}

\begin{definition} Let $T$ be a one-holed torus and let $\alpha$
be an interior simple closed geodesic of $T$. We denote
$h_\alpha$ the unique simple geodesic path which goes from
boundary to boundary and intersects boundary at two right angles
and does not cross $\alpha$. We will refer to the geodesic path
$h_\alpha$ as the height associated to $\alpha$ (see figure
\ref{fig:qpieceheight}).
\end{definition}

By using hyperbolic trigonometry, one can prove the following
result (for a proof, see for instance \cite{sc931}).

\begin{lemma}\label{lem:height} Let $T$ be a one-holed torus.
Let $\gamma$ be an interior simple closed geodesic of $T$ and
denote its associated height $h_\gamma$.  Then $\gamma$ passes
through exactly two of the three Weierstrass points and the
remaining Weierstrass point is the midpoint of $h_\gamma$.
Furthermore, the length of $\gamma$ is directly proportional to
the length of $h_\gamma$.
\end{lemma}

The following proposition, slightly more general than what we
require, has an interest in its own right.

\begin{proposition}\label{prop:torus2}
Let $T$ be a one-holed torus (where the boundary geodesic
$\varepsilon$ is allowed to be a cusp). Let $\alpha$ be an
interior simple closed geodesic and let $\beta$ be any other
interior simple closed geodesic that intersects $\alpha$ twice.
Then
\begin{equation*}
\ell(\beta)\geq 2 \,\arccosh \left( 1 +
\frac{\cosh\frac{\ell(\varepsilon)}{2}+1}{2\left(\cosh\frac{\ell(\alpha)}{2}-1\right)}\right).
\end{equation*}

Furthermore equality holds only when $T$ is obtained by pasting
$\alpha$ with a half-twist and $\beta$ is the shortest simple
closed geodesic that intersects $\alpha$ twice.
\end{proposition}
\begin{proof}
For a given $\alpha$ and $\varepsilon$, let $\beta$ be the
shortest simple closed geodesic $\beta$ that crosses $\alpha$
twice. Now consider the family of tori obtained by twisting along
$\alpha$. The key to the proof is showing that $\beta$ is shortest
when $\alpha$ is pasted with a half-twist. In first instance, let
us suppose that $\varepsilon$ is not a cusp. Consider the height
$h_\beta$ associated to $\beta$. By lemma \ref{lem:height}, the
length of $h_\beta$ is proportional to the length of $\beta$, so
minimizing the length of $\beta$ is equivalent to minimizing the
length of $h_\beta$. Denote by $e_1$, $e_2$ and $a$ the three
boundary to boundary geodesic perpendicular paths of the pair of
pants $(\alpha,\alpha,\varepsilon)$ as indicated on figure
\ref{fig:cuttorus}.

\begin{figure}[h]
\leavevmode \SetLabels
(.6*.8) $\alpha$\\
(.6*.1) $\alpha$\\
(.26*1.02) $\alpha$\\
(.25*.267) $a$\\
(.26*.03) $\varepsilon$\\
(.485*.45) $a$\\
(.88*.45) $a$\\
(.66*.45) $\varepsilon$\\
(.2*.21) $e_1$\\
(.295*.23) $e_2$\\
(.7*.63) $e_1$\\
(.7*.27) $e_2$\\
\endSetLabels
\begin{center}
\AffixLabels{\centerline{\epsfig{file=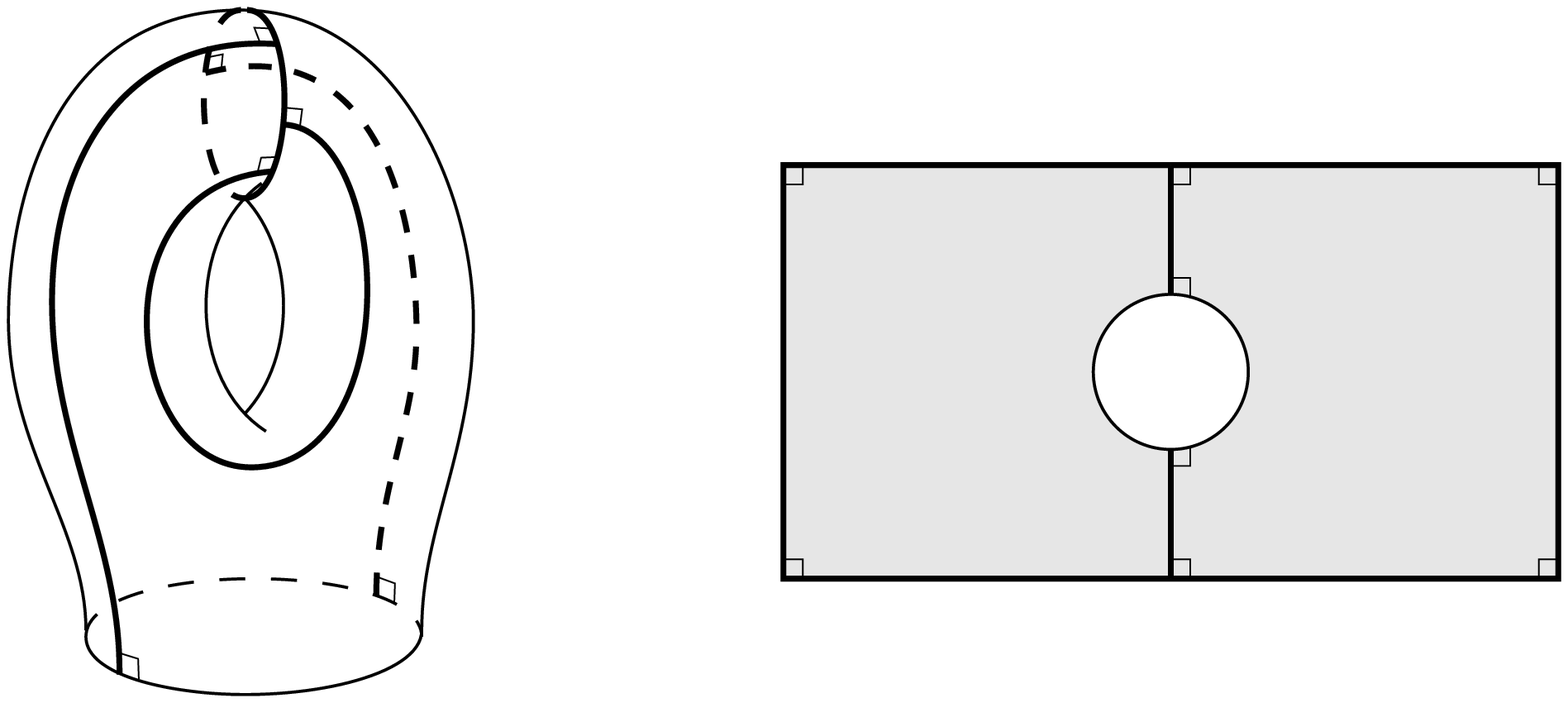,height=4.5cm}}}
\end{center}
\caption{The torus $T$ cut along $\alpha$ and $a$}
\label{fig:cuttorus}
\end{figure}

Cutting $T$ along $\alpha$ and path $a$ one obtains a one-holed
hyperbolic rectangle as in figure \ref{fig:cuttorus}. (This
particular way of viewing the one-holed torus is a central part of
\cite{busem88}.) Notice that $\ell(e_1)=\ell(e_2)$, which can be
seen either by using hyperbolic trigonometry or by using the
hyperelliptic involution. By cutting along paths $h_\alpha$,
$e_1$, $e_2$ and $a$, one would obtain four isometric right-angled
pentagons. The path $h_\beta$ intersects $\alpha$ twice, and thus
the two subpaths of $\beta$ between $\alpha$ and $\varepsilon$ are
of length at least $\ell(e_1)(=\ell(e_1))$, and the subpath from
$\alpha$ and back again is at least of length $\ell(a)$. Thus
$$\ell(h_\beta)\geq \ell(e_1)+\ell(e_2)+\ell(a).$$
Equality only holds when the path $h_\beta$ is exactly the path
$e_1\cup a \cup e_2$. This only occurs when the pasting is right,
meaning when $\alpha$ is pasted with a half-twist. Now, when
$\varepsilon$ is a cusp, we cannot immediately assume that the
optimal situation is when there is a half-twist, but this is true
because of the continuity of lengths of interior closed curves
when $\ell(\varepsilon)$ goes to $0$.\\

We now need to calculate the length of $\beta$ when $\alpha$ is
pasted with a half-twist. For this we shall use the well known
formulas for different types of hyperbolic polygons (see for
instance \cite[p.\;454]{bubook}). This can be done by considering
two hyperbolic polygons inscribed in $T$.

\begin{figure}[h]
\leavevmode \SetLabels
(.63*.785) $\alpha$\\
(.81*.8) {\scriptsize$\frac{\ell(\alpha)}{4}$}\\
(.82*.6) $Q$\\
(.565*.6) $\beta$\\
(.795*.26) $\beta$\\
(.6*.29) $P$\\
(.59*.07) {\scriptsize$\frac{\ell(\alpha)}{2}$}\\
(.73*.1) $\alpha$\\
(.28*.8) $\alpha$\\
(.22*.03) $\varepsilon$\\
(.475*.3) {\scriptsize$\frac{\ell(a)}{2}$}\\
(.885*.55) {\scriptsize$\frac{\ell(a)}{2}$}\\
(.76*.55) {\scriptsize$\frac{\ell(\beta)}{4}$}\\
(.67*.4) {\scriptsize$\frac{\ell(\varepsilon)}{4}$}\\
(.262*.1) $h_\alpha$\\
(.315*.45) $\beta$\\
\endSetLabels
\begin{center}
\AffixLabels{\centerline{\epsfig{file=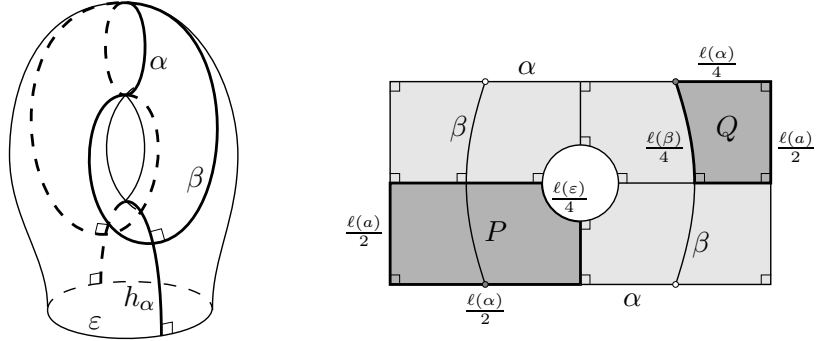,height=4.5cm}}}
\end{center}
\caption{The polygons $P$ and $Q$} \label{fig:torus2int}
\end{figure}

The first one, denoted $Q$, is one of the hyperbolic
quadrilaterals with three right angles delimited by arcs of paths
$a$, $\alpha$, $\beta$ and $h_\alpha$ as in figure
\ref{fig:torus2int}. The second polygon $P$ is one the four
isometric right-angled pentagons (or quadrilaterals with a point
at infinity when $\ell(\varepsilon)=0$) obtained by cutting $T$
along $\alpha$, $a$, $e_1$, $e_2$ and $h_\alpha$ (see figure
\ref{fig:torus2int}). Using the formulas for a quadrilateral with
three right angles, one obtains
$$
\sinh\frac{\ell(\beta)}{4} = \sinh \frac{\ell(a)}{2}\; \cosh
\frac{\ell(\alpha)}{4}.
$$
Now using the formula for a right-angled pentagon with $P$ we
obtain
$$
\sinh\frac{\ell(a)}{2}\; \sinh \frac{\ell(\alpha)}{2} = \cosh
\frac{\ell(\varepsilon)}{4}.
$$

Putting these two formulas together one obtains

$$
\sinh\frac{\ell(\beta)}{4}=\frac{\cosh\frac{\ell(\varepsilon)}{4}
\cosh\frac{\ell(\alpha)}{4}}{\sinh\frac{\ell(\alpha)}{2}}\;.
$$

With a little manipulation one obtains
$$
\cosh \frac{\ell(\beta)}{2}=1 +
\frac{\cosh\frac{\ell(\varepsilon)}{2}+1}{2\left(\cosh\frac{\ell(\alpha)}{2}-1\right)}\;,
$$

which proves the result.
\end{proof}

There is an immediate corollary which gives a universal lower
bound on the greatest of two lengths of two geodesics intersecting
twice on a one-holed torus.

\begin{corollary}\label{cor:torus2}
Let $T$ be a one-holed torus (where the boundary geodesic
$\varepsilon$ is allowed to be a cusp). Let $\alpha$ and $\beta$
be two interior simple closed geodesics that intersect twice. Then
$$
\max\{\ell(\alpha),\ell(\beta)\}\geq 2\,\arccosh 2.
$$

Furthermore equality holds for a torus $T$ with a cusp which
contains a simple closed geodesic $\alpha$ of length $2\,\arccosh
2$ pasted with a half-twist, and taking $\beta$ to be the shortest
simple closed geodesic which intersects $\alpha$ twice.
\end{corollary}

\begin{proof}
Suppose that $\ell(\alpha) < 2\,\arccosh 2$. By proposition
\ref{prop:torus2}, $\ell(\beta)>2\,\arccosh 2$. Now if
$\ell(\alpha) = 2\,\arccosh 2$, by proposition \ref{prop:torus2},
$\ell(\beta)=2\,\arccosh 2$ if and only if $\ell(\varepsilon)=0$
and $T$ is obtained by pasting $\alpha$ with a half-twist.
\end{proof}

Note that the torus described in corollary \ref{cor:torus2} is the
same torus as $S_1$. To see this, we shall find two simple closed
geodesics that intersect once, and both of length
$2\,\arccosh\sqrt{2}$. Consider the quadrilateral $Q$ in figure
\ref{fig:torus2int} and in particular the diagonal of $Q$ from top
left to bottom right. Now consider the diagonals of each one of
the four isometric copies of $Q$. Together these four geodesic
paths form two simple closed geodesics, say $\gamma_1$ and
$\gamma_2$, of equal length that intersect once. When
$\ell(\epsilon)=0$ and $\ell(\alpha)=2\,\arccosh 2$, a quick
calculation shows that
$\ell(\gamma_1)=\ell(\gamma_2)=2\,\arccosh\sqrt{2}$. As $S_1$ is
unique up to isometry, the two tori are the same.\\

\begin{theorem}\label{thm:S_2}
The surface $S_2$ is the one-holed torus described in corollary
\ref{cor:torus2} and $L_2 = 2\,\arccosh(2)$.
\end{theorem}

\begin{proof}
In view of proposition \ref{prop:fourholed2} and corollary
\ref{cor:torus2}, we now know that $S_2$ is a torus with one or
two punctures. Suppose $S_2$ is a torus with two punctures, i.e.,
the curves labeled $\tilde{\varepsilon}_1$ and
$\tilde{\varepsilon}_2$ on figure \ref{fig:AuBu} are homotopic to
cusps.

\begin{figure}[ht!]
\centerline{\SetLabels
(.25*.45) {$\delta$}\\
(.765*.45) {$\delta$}\\
(.5*.99) {$\beta$}\\
(.16*.015) {$x$}\\
(.16*.935) {$y$}\\
(.4*.015) {$x'$}\\
(.4*.935) {$y'$}\\
(.61*.015) {$y$}\\
(.61*.935) {$x$}\\
(.85*.015) {$y'$}\\
(.85*.935) {$x'$}\\
(.5*-0.05) {$\beta$}\\
(.64*.26) {$\gamma$}\\
(.55*.21) {$a_1$}\\
(.01*.45) {$a_2$}\\
(1.01*.45) {$a_2$}\\
\endSetLabels
\AffixLabels{
\hspace{2mm}\psfig{figure=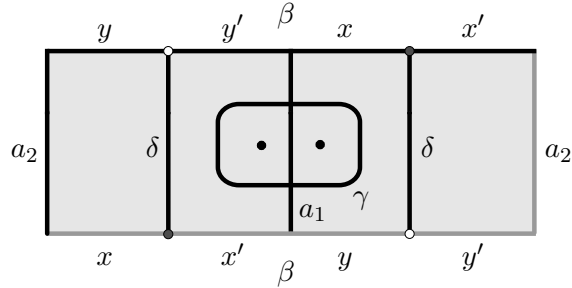,height=3.3cm}\hspace{2mm}}}
\caption{Two intersections on a twice-punctured torus}
\label{fig:fish2int}
\end{figure}

Now suppose we have two simple closed geodesics, $\alpha$ and
$\beta$, with $\ell(\alpha)\geq \ell(\beta)$ and
$\ell(\alpha)\leq 2 \,\arccosh 2$. (If this is not possible, then
necessarily $S_2$ is the once punctured torus of corollary
\ref{cor:torus2}.) The geodesic $\alpha$ is cut into two arcs by
$\beta$, say $a_1$ and $a_2$. Suppose $\ell(a_1)\geq \ell(a_2)$.
Consider the geodesic curves $\gamma$ and $\delta$ as in figure
\ref{fig:fish2int}. $\gamma$ is the separating curve that
intersects $a_1$ twice but doesn't intersect $a_2$ or $\beta$,
and $\delta$ is the curve that intersects $\beta$ twice but
doesn't intersect $\gamma$ or $\alpha$. Consider the lengths
$x,x',y,y'$ of the different arcs of $\beta$ as labeled on figure
\ref{fig:fish2int}. We have $\ell(\delta)< x + \ell(a_2) + y + x'
+ \ell(a_2) +y'= \ell(\alpha)+\ell(\beta)\leq 2 \ell(\alpha)$.
Notice that this implies the width of the collar around $\delta$
satisfies
$$
w(\delta)>\arcsinh\left(\frac{1}{\sinh \ell(\alpha)}\right). $$

We can now apply the collar theorem to $\beta$, using the fact
that $\beta$ intersects both $\alpha$ and $\delta$ twice and
$\alpha$ and $\delta$ do not intersect. The collar theorem
\ref{thm:collar} implies that the length of $\beta$ satisfies the
following inequality:
$$
\ell(\beta)\geq 4 w(\alpha) + 4 w(\delta) > 2 \,\arccosh 2\geq
\ell(\alpha).
$$

This proves the result.
\end{proof}

\section{Finding $S_3$ and calculating $L_3$}

Let $\alpha$ and $\beta$ be two simple closed geodesics on a
Riemann surface that intersect three times. Name the intersection
points $A$, $B$ and $C$ and orient $\alpha$ and $\beta$ such that
$A$, $B$ and $C$ come in that order on $\alpha$ and on $\beta$. As
in the case of two intersections, we consider a collar around
$\beta$ and the orientations induced on it in the different
intersection points by the ordered pair of simple closed
\emph{oriented} geodesics $(\alpha,\beta)$. We distinguish two
situations:
\begin{enumerate}
\item $(\alpha,\beta)$ induces opposite orientations in two of the three
intersection points (without loss of generality we can assume that
$(\alpha,\beta)$ induces opposite orientations in $A$ and in $B$),

\item $(\alpha,\beta)$ induces the same orientation in $A$, in $B$ and in $C$.
\end{enumerate}

In the first situation, lemma \ref{lem:3pseudoX} will show that
$\max\{\ell(\alpha),\ell(\beta)\}\geq2\,\arccosh(3)$.\\

In the second situation, we will show that the optimal surface is
a torus with a cusp containing two simple closed geodesics of
lengths approximately $2\,\arccosh(2.648)$ intersecting one
another three times.

\begin{lemma}\label{lem:3pseudoX}
Let $\alpha$ and $\beta$ be two simple closed oriented geodesics
on a Riemann surface that intersect three times in $A$, $B$ and
$C$, such that $A,B,C$ are consecutive on both $\alpha$ and
$\beta$.

If the ordered pair $(\alpha,\beta)$ induces opposite orientations
on the surface in $A$ and in $B$, then
$\max\{\ell(\alpha),\ell(\beta)\}\geq2\,\arccosh(3)$.
\end{lemma}

\begin{proof}
Comparing the lengths of the arcs between $B$ and $C$, there are
two possible situations:
\begin{enumerate}
\item The length $\overline{BC}_\alpha$ of the oriented geodesic arc from
$B$ to $C$ on the geodesic $\alpha$ is smaller then
$\overline{BC}_\beta$, the length of the oriented geodesic arc
from $B$ to $C$ on the geodesic $\beta$

\item
This is not the case, meaning $\overline{BC}_\alpha\geq
\overline{BC}_\beta$.
\end{enumerate}

We now build the oriented closed curves $\tilde\gamma$ and
$\tilde\delta$:
\begin{itemize}
\item In situation 1, we set $\tilde\gamma=\alpha$; in situation
2, $\tilde\gamma$ is obtained following $\alpha$ from $A$ to $B$,
then $\beta$ from $B$ to $C$ and again $\alpha$ from $C$ to $A$.
\item In situation 1, $\tilde\delta$ is obtained following $\beta$
from $A$ to $B$, then $\alpha$ from $B$ to $C$ and again $\beta$
from $C$ to $A$; in situation 2, we set $\tilde\delta=\beta$.
\end{itemize}

These two curves $\tilde\gamma$ and $\tilde\delta$ are thus
homotopic to two simple closed oriented geodesics $\gamma$ and
$\delta$ intersecting one another twice such that
$$
\max\{\ell(\alpha),\ell(\beta)\}\geq\max\{\ell(\gamma),\ell(\delta)\}.
$$

Furthermore, the ordered pair $(\gamma,\delta)$ induces opposite
orientations in its two intersection points. Therefore
$\max\{\ell(\gamma),\ell(\delta)\}\geq 2\,\arccosh(3)$ by
proposition \ref{prop:fourholed2}.
\end{proof}

\begin{lemma}\label{l:reduce-to-Q}
Let $S$ be a Riemann surface and let $\alpha$ and $\beta$ be two
oriented simple closed geodesics on $S$ intersecting one another
three times such that the ordered pair $(\alpha,\beta)$ induces
the same orientation on $S$ in every intersection. Name the
intersections $A,B,C$ such that they are consecutive on $\alpha$.
If $A,B,C$ are also consecutive on $\beta$, then there is a torus
with one cusp or a torus with two cusps containing two simple
closed geodesics $\gamma$ and $\delta$ which satisfy
$\ii(\alpha,\beta)=3$ and
$\max\{\ell(\alpha),\ell(\beta)\}\geq\max\{\ell(\gamma),\ell(\delta)\}$.
\end{lemma}

\begin{proof}
The surface $S_{\alpha,\beta}$ is of signature $(1,3)$, $(1,2)$ or
$(1,1)$, as can easily be seen by considering whether the simple
closed curves $\tilde\varepsilon_1$, $\tilde\varepsilon_2$ and
$\tilde\varepsilon_3$ (shown in figure \ref{fig:AuBuCu}) are
null-homotopic or not.

\begin{figure}[ht!]
\centerline{\SetLabels
(.19*.45) {?}\\
(.5*.45) {?}\\
(.81*.45) {?}\\
(.02*.11) {$\alpha$}\\
(.375*.11) {$\alpha$}\\
(.68*.11) {$\alpha$}\\
(.2*.99) {$\beta$}\\
(.05*-.08) {$A$}\\
(.05*.99) {$B$}\\
(.95*-.08) {$A$}\\
(.95*.99) {$B$}\\
(.35*.99) {$C$}\\
(.35*-.08) {$B$}\\
(.65*.99) {$A$}\\
(.65*-.08) {$C$}\\
(.11*.27) {$\tilde\varepsilon_1$}\\
(.415*.27) {$\tilde\varepsilon_2$}\\
(.721*.27) {$\tilde\varepsilon_3$}\\
\endSetLabels
\AffixLabels{
\hspace{2mm}\psfig{figure=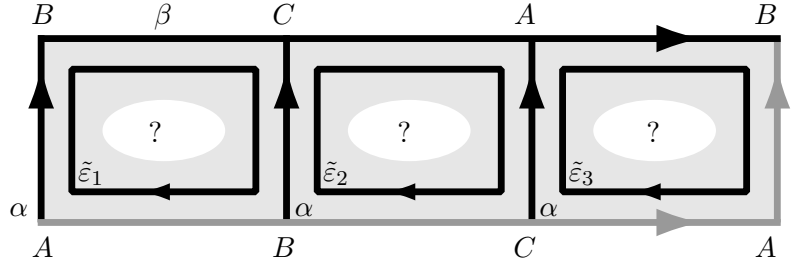,height=2.9cm}\hspace{2mm}}}
\caption{The simple closed curves $\tilde\varepsilon_1$,
$\tilde\varepsilon_2$ and $\tilde\varepsilon_3$}\label{fig:AuBuCu}
\end{figure}

If one of the curves $\tilde\varepsilon_1$, $\tilde\varepsilon_2$
or $\tilde\varepsilon_3$ is null-homotopic, corollary
\ref{cor:cusp} proves the lemma. Otherwise, the optimal
topological situation is a torus with three cusps (again due to
corollary \ref{cor:cusp}). On this surface, there is a simple
closed geodesic $\eta$ dividing the surface into $X_\eta$, a
sphere with three cusps and  boundary geodesic $\eta$, and
$T_\eta$, a surface of signature $(1,1)$. Notice that $\beta$ is
entirely contained in $T_\eta$ as can be seen in figure
\ref{fig:XQ}.\enlargethispage{3mm}

\begin{figure}[h]
\centerline{ \SetLabels
(.0*.7) {$A$}\\
(.33*.7) {$B$}\\
(.31*.8) {$U$}\\
(.31*.9) {$V$}\\
(.66*.7) {$C$}\\
(.69*.8) {$Y$}\\
(.705*.9) {$W$}\\
(.99*.7) {$A$}\\
(.99*1) {$B$}\\
(-.01*.85) {$\alpha$}\\
(.0*1) {$B$}\\
(.33*1) {$C$}\\
(.66*1) {$A$}\\
(.45*1) {$\beta$}\\
(.56*.8) {$\eta$}\\
(.99*.31) {$\beta$}\\
(.22*.155) {$\alpha$}\\
(.2*.26) {$\varepsilon$}\\
(.51*.415) {$\eta$}\\
(.63*.39) {$Y$}\\
(.43*.15) {$V$}\\
(.62*.14) {$U$}\\
(.44*.47) {$W$}\\
\endSetLabels
\AffixLabels{\vspace{1mm}\psfig{figure=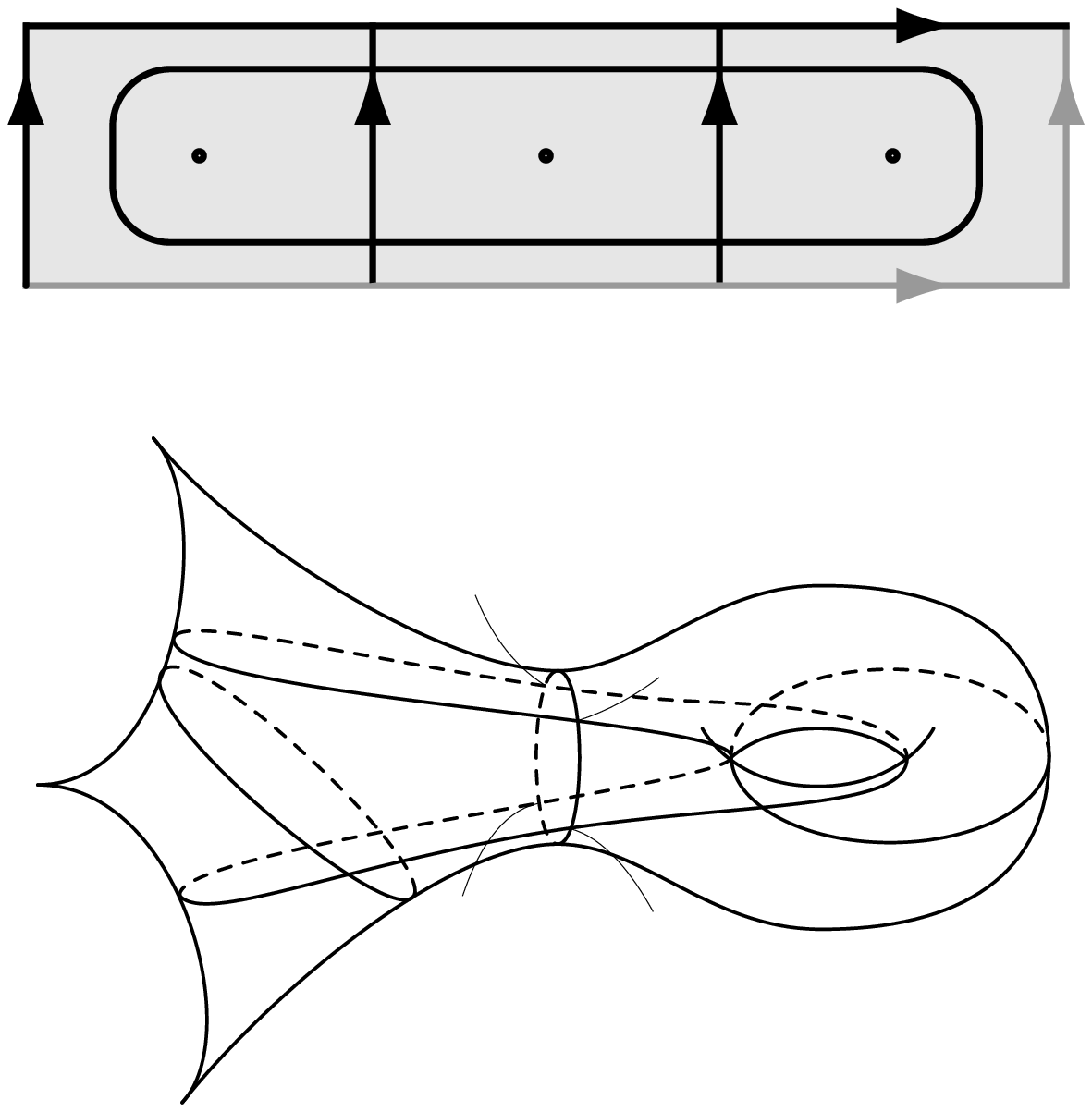,height=9cm}}}
\caption{Three intersections on a torus with three
cusps}\label{fig:XQ}
\end{figure}

The intersection points between $\alpha$ and $\eta$ will be
denoted $U$, $V$, $W$ and $Y$ as in figure \ref{fig:XQ}. First
consider the geodesic arc of $\alpha$ from $Y$ to $W$. There is a
dividing geodesic $\varepsilon$ on $X_\eta$, that does not
intersect this arc. Cutting $X_\eta$ along $\varepsilon$, we get a
surface of signature $(0,3)$. We can now diminish the length of
$\varepsilon$ in order to get another cusp. This surface of
signature $(0,3)$ with two cusps and the boundary geodesic $\eta$
contains a geodesic arc from $Y$ to $W$ that is shorter than the
original arc from $Y$ to $W$ on $X_\eta$ (this is part of the
statement of the technical lemma used in \cite{pa044} in order to
show theorem \ref{thm:cusps}).

Obviously, we can do the same for the geodesic arc joining $U$ and
$V$. Thus we can replace $X_\eta$ by the unique surface of
signature $(0,3)$ with two cusps and a boundary geodesic of length
$\ell(\eta)$. We get a torus with two cusps that contains a
geodesic $\beta$ and a curve $\tilde \alpha$ that intersect three
times and such that $\ell(\alpha)\geq \ell(\tilde \alpha)$.
Therefore, the geodesic $\gamma$ that is homotopic to $\tilde
\alpha$ intersects the geodesic $\beta$ (that we rename to
$\delta$) three times and
$\max\{\ell(\alpha),\ell(\beta)\}\geq\max\{\ell(\gamma),\ell(\delta)\}$.
\end{proof}

\begin{lemma}
\label{l:three-intersections-on-Q} If $\alpha$ and $\beta$ are two
simple closed geodesics on a one-holed torus satisfying
$\ii(\alpha,\beta)=3$, then
$$
\textstyle\max\{\ell(\alpha),\ell(\beta)\}\geq2\,\arccosh
\left(\sqrt{\frac12\left(7+\frac{11}{3}\sqrt{\frac{11}{3}}\right)}\;\right).
$$

This bound is sharp and is reached by a unique once-punctured
torus up to isometry.
\end{lemma}
\begin{proof}
We shall use the parameters for the set of isometry classes of
one-holed tori found in \cite{busem88}. Let $(r,s,t)$ be a set of
these parameters such that $1<r\leq s\leq t\leq rs$ where $r$, $s$
and $t$ are the half-traces (hyperbolic cosines of half of the
lengths) of the shortest three geodesics $\varrho$, $\sigma$ and
$\tau=(\varrho\sigma)^{-1}$. (In \cite{busem88}, half-traces are
denoted traces, but we shall continue to use the term half-traces
as it is more standard.) Then, the geodesics
$\alpha=\varrho\sigma^{-1}$ and $\beta=\tau\varrho^{-1}$ intersect
three times and $\alpha$ is the forth shortest simple closed
geodesic (see \cite{busem88} for details). The half-traces of
$\alpha$ and $\beta$ are $a=2rs-t$ and $b=2rt-s$.

For a fixed $r$, $\max\{a,b\}=b=2rt-s$ is therefore minimal if
$s=t$. In this case $0=2rst-r^2-s^2-t^2=2s^2(r-1)-r^2$ and
therefore $b^2=s^2(2r-1)^2=\frac{r^2(2r-1)^2}{2(r-1)}$.

But for $r>1$, this last quantity is minimal for
$$\textstyle\frac{d}{dr}\frac{r^2(2r-1)^2}{2(r-1)}=0\ \Longleftrightarrow\
\frac{r(2r-1)(6r^2-9r+2)}{2(r-1)^2}=0,$$ i.e.,
$r=\frac14\left(3+\sqrt{\frac{11}{3}}\right)$. Therefore
$b^2=s^2(2r-1)^2
\geq\frac12\left(7+\frac{11}{3}\sqrt{\frac{11}{3}}\right)$.\\

There is a a torus with one cusp
on which there are two geodesics of lengths $2\,\arccosh
\left(\sqrt{\frac12\left(7+\frac{11}{3}\sqrt{\frac{11}{3}}\right)}\;\right)$
intersecting one another three times. Its trace coordinates are
necessarily $(r,s,t) =
\left(\frac14\left(3+\sqrt{\frac{11}{3}}\right),
\sqrt{\frac{13+7\sqrt{\frac{11}{3}}}{8}},
\sqrt{\frac{13+7\sqrt{\frac{11}{3}}}{8}}\right)$, up to a choice
of curves with trace $(r,s,t)$. Therefore the bound is sharp and
is attained by a unique once-punctured torus up to isometry.
\end{proof}

It is worth noticing the the torus described in this lemma is {\it
not} $S_1$. As mentioned in the proof, its systole length is $2\,
\arccosh\left(\frac14\left(3+\sqrt{\frac{11}{3}}\right)\right)$
and not $2\,\arccosh\sqrt{2}$.\\

\begin{theorem}\label{thm:S_3}
The surface $S_3$ is the one-holed torus described in lemma
\ref{l:three-intersections-on-Q}  and $$\textstyle L_3 =
2\,\arccosh
\left(\sqrt{\frac12\left(7+\frac{11}{3}\sqrt{\frac{11}{3}}\right)}\;\right).$$
\end{theorem}

\begin{proof}
By what precedes, $S_3$ is a torus with one or two punctures. As
in the proof of theorem \ref{thm:S_2}, let us suppose that there
exists a twice-punctured torus $T$ with two geodesics $\alpha$ and
$\beta$ that intersect three times, and both of length less or
equal to $2\,\arccosh
\left(\sqrt{\frac12\left(7+\frac{11}{3}\sqrt{\frac{11}{3}}\right)}\;\right)$.
For the remainder of the proof, denote this constant $k_3$, as
$L_3$ would be premature.

\begin{figure}[ht!]
\centerline{\SetLabels
(.2*.99) {$b_1$}\\
(.2*-.08) {$b_3$}\\
(.51*.99) {$b_2$}\\
(.51*-.08) {$b_1$}\\
(.82*.99) {$b_3$}\\
(.82*-.08) {$b_2$}\\
(.01*.45) {$a_3$}\\
(.32*.45) {$a_1$}\\
(.7*.45) {$a_2$}\\
(1.01*.45) {$a_3$}\\
(.4*.71) {$\gamma_\beta$}\\
(.6*.45) {$\gamma_\alpha$}\\
(.1*.45) {$\gamma_\alpha$}\\
\endSetLabels
\AffixLabels{
\hspace{2mm}\psfig{figure=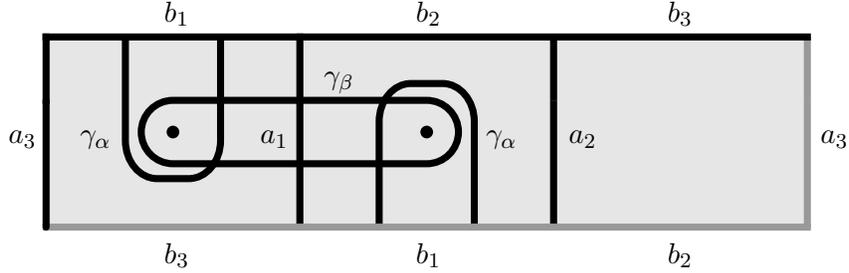,height=3cm}\hspace{2mm}}}
\caption{Three intersections on a twice-punctured torus
}\label{fig:fish3inta}
\end{figure}


Both $\alpha$ and $\beta$ are separated into three paths by each
other, and let us denote these paths respectively $a_1$, $a_2$ and
$a_3$ for $\alpha$ and $b_1$, $b_2$ and $b_3$ for $\beta$.
The pasting condition implies that we are now in the situation
illustrated in figure \ref{fig:fish3inta}. On this figure, two
additional simple closed curves have been added, and are denoted
$\gamma_\alpha$ and $\gamma_\beta$. The curve $\gamma_\alpha$ is
defined as the unique separating simple closed geodesic that does
not intersect $\alpha$ and intersects $\beta$ minimally (twice),
and $\gamma_\beta$ is defined symmetrically. We will use a rough
upper-bound on the sum of their lengths. It is easy to see that
$$
\ell(\gamma_\alpha)+\ell(\gamma_\beta)<
(2\ell(a_1)+\ell(a_2)+\ell(a_3)+\ell(b_2) +\ell(b_3))
$$
$$
+(\ell(a_2)+\ell(a_3)+2\ell(b_1)+\ell(b_2)+\ell(b_3))
$$
$$
= 2\ell(\alpha)+2\ell(\beta).
$$

This implies that
$\min\{\ell(\gamma_\alpha),\ell(\gamma_\beta\})\leq
2\max\{\alpha,\beta\}\leq 2k_3$. So far, we have made no
particular assumptions on $\alpha$ and $\beta$, so without loss
of generality we can suppose that $\alpha$ is such that
$\ell(\gamma_\alpha)\leq 2k_3$.

\begin{figure}[ht!]
\centerline{\SetLabels
(.25*.535) {$\alpha$}\\
(.95*.495) {$\alpha'$}\\
(.495*.5) {$\gamma_\alpha$}\\
(.74*.725) {$\beta$}\\
\endSetLabels
\AffixLabels{\psfig{figure=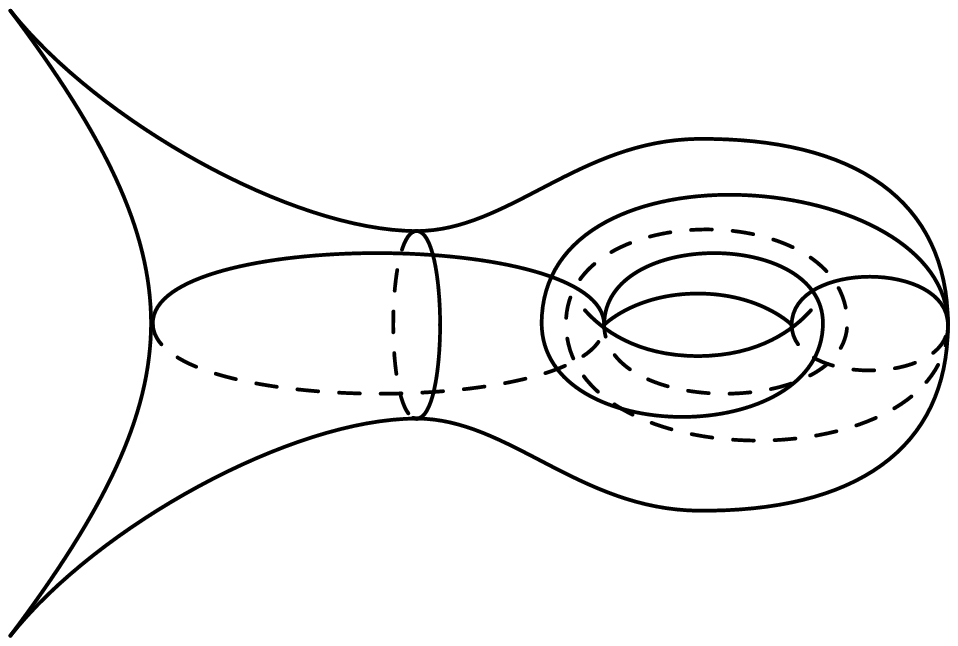,height=5.5cm}}}
\caption{}\label{fig:fish3intb}
\end{figure}

Denote by $\alpha'$ the unique simple closed geodesic of $T$ that
intersects neither $\alpha$ nor $\gamma_\alpha$. Notice that
$\alpha'$ intersects $\beta$ three times. We shall now find an
upper-bound on the length of $\alpha'$. By cutting along $\alpha$
and $\alpha'$, one obtains two (isometric) pairs of pants.
Consider one of them as in figure \ref{fig:fish3intc}. We denote
by $h_\alpha$ the shortest non-trivial path from $\alpha$ and
back again. Notice that
\begin{equation}\label{eqn:k3}
\ell(h_\alpha)\leq \frac{\ell(\gamma_\alpha)}{2}\leq k_3.
\end{equation}

\begin{figure}[ht!]
\centerline{\SetLabels
(.32*.53) {$h_\alpha$}\\
(.44*.82) {$\alpha'$}\\
(1*.72) {$\frac{\ell(\alpha')}{2}$}\\
(.2*.01) {$\alpha$}\\
(.75*.02) {$l$}\\
(.88*.02) {$l'$}\\
(.81*-.06) {$\frac{\ell(\alpha)}{2}$}\\
(.845*.38) {$\frac{\ell(h_\alpha)}{2}$}\\
\endSetLabels
\AffixLabels{\psfig{figure=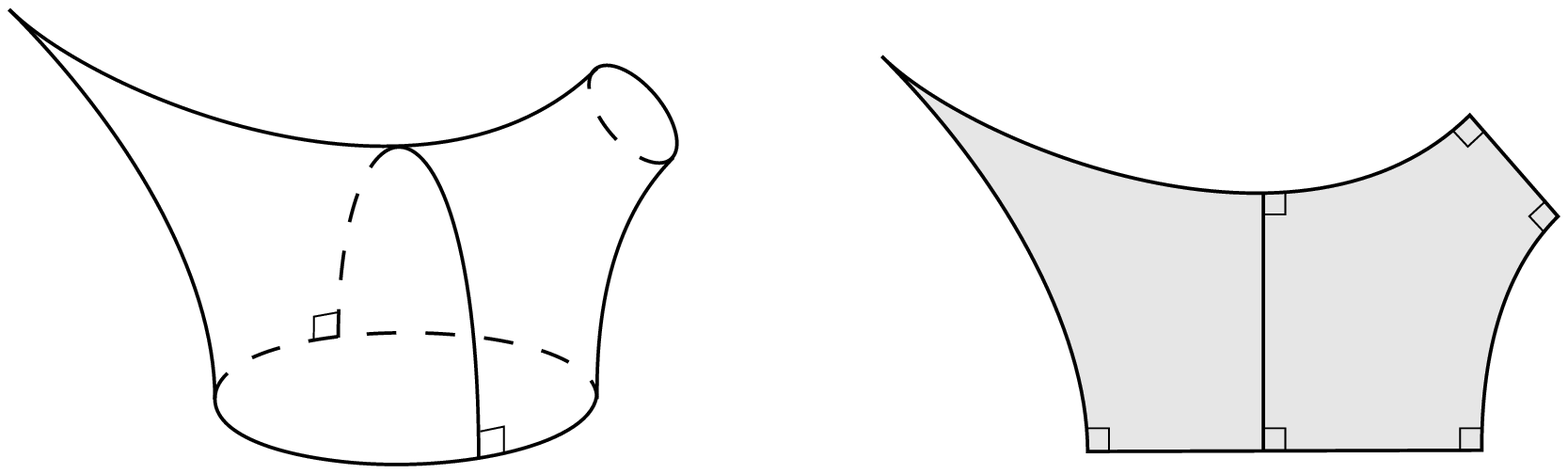,height=4cm}}}
\caption{Bounding the length of $\alpha'$}\label{fig:fish3intc}
\end{figure}

Consider the lengths $l$ and $l'$ in figure \ref{fig:fish3intc}.
Once again, we shall make use of the formulas for hyperbolic
polygons. Using the hyperbolic trigonometry formulas for a
pentagon with right angles, we obtain
$$
\sinh \frac{\ell(h_\alpha)}{2}\; \sinh l'= \cosh
\frac{\ell(\alpha')}{2}.
$$

Using the formulas for a quadrilateral with three right angles and
one zero angle, one obtains $ \sinh l\; \sinh
\frac{\ell(h_\alpha)}{2} =1$ and  equation \ref{eqn:k3} leads to
$$
l > \frac{1}{\sinh \frac{k_3}{2}}\;.
$$

From these equations, and our initial hypothesis on the lengths of
$\alpha$ and $\beta$, we obtain the following bound on the length
of $\alpha'$:
$$
\ell(\alpha') < 2 \,\arccosh\left( \sinh \frac{k_3}{2}\;
\sinh\left(\frac{k_3}{2}- \arcsinh \frac{1}{\sinh
\frac{k_3}{2}}\right)\right).
$$

This implies that the collar width of $\alpha'$ satisfies
$w(\alpha')>0.25$. As $\beta$ intersects both $\alpha$ and
$\alpha'$ three times and $\alpha$ and $\alpha'$ are disjoint, we
have that
$$
\ell(\beta)\geq 6 w(\alpha) + 6 w(\alpha').
$$

As $\ell(\alpha)\leq k_3$ implies $w(\alpha)>0.3$, we now have
$\ell(\beta) > 1.5 + 1.8 > k_3$ which contradicts the hypotheses.
Thus $S_3$ is a once-punctured torus and we can apply lemma
\ref{l:three-intersections-on-Q}.
\end{proof}

\section{Concluding remarks}

The surfaces $S_1=S_2$ and $S_3$ are specific once-punctured tori.
Both admit automorphisms distinct from the hyperelliptic
involution. $S_1$ admits a number of automorphisms both conformal
and anticonformal. Using the main result of \cite{copa051}, $S_3$
admits an orientation reversing involution because it can be
obtained by pasting a simple closed geodesic with a half-twist,
but does not admit a non-trivial conformal automorphism. This is
not so surprising seeing as there are only two isometry classes of
once-punctured tori that admit a non-trivial conformal
automorphism, namely $S_1$ and the torus
with largest automorphism group, often called the Modular torus.\\

Finding $S_k$ for $k\geq 4$ seems like a difficult problem, but
can we say something about the set of $S_k$? For higher
intersection number, it is not clear whether or not $S_k$ even has
boundary (recall that two simple closed geodesics can fill closed
surfaces if they are allowed sufficiently many intersection
points). In spite of this remark, it seems reasonable to
conjecture that $S_k$ is always a once-punctured torus.
Furthermore, due to the existence on $S_k$ of geodesics of equal
length, it also seems reasonable to conjecture that the $S_k$ all
have non-trivial automorphism groups. Supposing that the $S_k$ are
all once-punctured tori, are they all found in a finite set of
isometry classes of once-punctured tori?\\

\noindent{\bf Acknowledgement.}\\

\noindent The authors would like to thank Ying Zhang for his
comments.

\bibliographystyle{plain}
\def\cprime{$'$}

\end{document}